\documentclass[a4paper,12pt]{amsart}
\usepackage{amssymb,amscd,amsmath,a4wide,graphicx,stmaryrd,fullpage,setspace,microtype}
\usepackage[pagebackref=true]{hyperref}
\hypersetup{backref}
\usepackage[all]{xy}
\title{Left determined model categories}

\author[P. Gaucher]{Philippe Gaucher}

\address{Laboratoire PPS  (CNRS UMR 7126)\\
  Univ Paris Diderot\\
Sorbonne Paris Cit\'e\\
 Case 7014\\ 75205 PARIS Cedex 13 \\ France}

\urladdr{http://www.pps.univ-paris-diderot.fr/{\~{}}gaucher/} 



\swapnumbers

\newcommand{\C}{\mathcal{C}}

\newcommand{\K}{\mathcal{K}}
\newcommand{\LL}{\mathcal{A}}

\newcommand{\W}{\mathcal{W}}
\newcommand{\F}{\mathcal{F}}

\newcommand{\p}\times

\newtheorem{thm}{Theorem}[section]
\newtheorem{prop}[thm]{Proposition}
\newtheorem{lem}[thm]{Lemma}
\newtheorem{cor}[thm]{Corollary}

\newtheorem{defn}[thm]{Definition}
\newtheorem{nota}[thm]{Notation}

\newcommand{\bd}{\begin{defn}}
\newcommand{\ed}{\end{defn}}
\newcommand{\bp}{\begin{prop}}
\newcommand{\ep}{\end{prop}}
\newcommand{\bth}{\begin{thm}}
\renewcommand{\eth}{\end{thm}}
\newcommand{\bpf}{\begin{proof}}
\newcommand{\epf}{\end{proof}}
\newcommand{\bc}{\begin{cor}}
\newcommand{\ec}{\end{cor}}

\newcommand{\fL}[1]{\ar@{->}[ll]_-{#1}}
\newcommand{\fR}[1]{\ar@{->}[rr]^-{#1}}
\newcommand{\fRr}[1]{\ar@{->}[rrr]^-{#1}}
\newcommand{\fD}[1]{\ar@{->}[dd]_-{#1}}
\newcommand{\fU}[1]{\ar@{->}[uu]^-{#1}}
\newcommand{\f}[2]{\ar@{->}[#1]|{#2}}
\newcommand{\ff}[2]{\ar@2{->}[#1]|{#2}}
\newcommand{\frr}[1]{\ar@{->}[rrrr]^-{#1}}

\newcommand{\fl}[1]{\ar@{->}[l]_-{#1}}
\newcommand{\fr}[1]{\ar@{->}[r]^-{#1}}
\newcommand{\fd}[1]{\ar@{->}[d]_-{#1}}
\newcommand{\fu}[1]{\ar@{->}[u]^-{#1}}

\newcommand{\iso}{\cong}

\renewcommand{\leq}{\leqslant}
\renewcommand{\geq}{\geqslant}

\def\cartesien{%
  \ar@{-}[]+R+<6pt,-2pt>;[]+RD+<6pt,-6pt>%
  \ar@{-}[]+D+<2pt,-6pt>;[]+RD+<6pt,-6pt>%
}
\def\cocartesien{%
  \ar@{-}[]+L+<-6pt,+2pt>;[]+LU+<-6pt,+6pt>%
  \ar@{-}[]+U+<-2pt,+6pt>;[]+LU+<-6pt,+6pt>%
}
\def\hocartesien{%
  \ar@{-}[]+R+<6pt,-2pt>;[]+RD+<6pt,-6pt>_{h}%
  \ar@{-}[]+D+<2pt,-6pt>;[]+RD+<6pt,-6pt>%
}
\def\hococartesien{%
  \ar@{-}[]+L+<-6pt,+2pt>;[]+LU+<-6pt,+6pt>_{h}%
  \ar@{-}[]+U+<-2pt,+6pt>;[]+LU+<-6pt,+6pt>%
}

\newcommand{\brm}[1]{\rm{\mathbf{#1}}}

\DeclareMathOperator{\id}{Id}

\DeclareMathOperator{\Mor}{Mor}

\newcommand{\wts}{\mathcal{W\!T\!S}}

\DeclareMathOperator{\bl}{\brm{\underline{L}}}

\makeatletter
\def\varholim@#1#2{%
  \vtop{\m@th\ialign{##\cr
    \hfil$#1\operator@font holim$\hfil\cr
    \noalign{\nointerlineskip\kern1.5\ex@}#2\cr
    \noalign{\nointerlineskip\kern-\ex@}\cr}}%
}
\def\holimproj{%
  \mathop{\mathpalette\varholim@{\leftarrowfill@\textstyle}}\nmlimits@
}
\def\holiminj{%
  \mathop{\mathpalette\varholim@{\rightarrowfill@\textstyle}}\nmlimits@
}
\makeatother

\DeclareMathOperator{\cof}{{\brm{cof}}}

\DeclareMathOperator{\inj}{{\brm{inj}}}
\newcommand{\ddownarrow}{{\downarrow}}

\DeclareMathOperator{\cyl}{{Cyl}}
\DeclareMathOperator{\cocyl}{{Path}}
\begin{document}

\begin{abstract}
  Several methods for constructing left determined model structures
  are expounded.  The starting point is Olschok's work on locally
  presentable categories. We give sufficient conditions to obtain left
  determined model structures on a full reflective subcategory, on a
  full coreflective subcategory and on a comma category.  An
  application is given by constructing a left determined model
  structure on star-shaped weak transition systems.
\end{abstract}

\maketitle

\tableofcontents

\section{Introduction}

\subsection*{Summary}

The notion of combinatorial model category is a powerful framework for
doing homotopy \cite{MR1780498} \cite{MR2506258}. It consists of a
locally presentable category equipped with a cofibrantly generated
model structure. Among them, there are the left determined ones in the
sense of \cite{rotho}, that is the combinatorial model categories such
that the class of weak equivalences is minimal with respect to a given
class of cofibrations. The interest of constructing left determined
model structures is that, for a given class of cofibrations, all other
ones are left Bousfield localizations of the left determined
one. J. H. Smith conjectured that for any locally presentable category
and any set of maps $I$, there exists a left determined combinatorial
model category such that the class of cofibrations is generated by
$I$. This statement is a consequence of Vop\v{e}nka's principle
\cite[Theorem~2.2]{rotho}. To our knowledge, the conjecture is still
open without assuming this large-cardinal axiom.

A remarkable step towards a proof of this conjecture is Olschok's
paper \cite{MO}. The latter paper generalizes Cisinski's work about
the homotopy theory of toposes \cite{MR1924082} to the framework of
locally presentable categories. It proves the existence of this left
determined model structure under reasonable hypotheses.

Several model structures are constructed in \cite{cubicalhdts},
\cite{homotopyprecubical} and \cite{csts} by using Olschok's work. The
common pattern of all these constructions is to start from an
application of Oslchok's theorem and to restrict the model structure
to reflective and coreflective full subcategories.

We expound here in full generality these methods. This paper is
written for two reasons: 1) we will use these methods repeatedly in
our studies of higher dimensional transition systems~\footnote{This is
  a work in progress belonging to the interface between algebraic
  topology and concurrency theory in computer science}, in particular
in the companion paper \cite{biscsts1}, 2) we hope that some people
will find these methods useful and maybe generalizable. This paper is
therefore designed to be a toolbox. Not only are methods for obtaining
left determined model structures on reflective and coreflective
subcategories given in this paper, but also sufficient conditions for
the standard model structure on a comma category to be left determined
as well. This paper ends with an application to star-shaped weak
transition systems.

\subsection*{Outline of the paper}

Section~\ref{good-model-cat} recalls Olschok's work and introduces the
notion of Olschok model structure. When the associated cartesian
cylinder is very good, we obtain a left determined model structure by
choosing an empty set of generating anodyne cofibrations. There is
nothing new in this section except Proposition~\ref{model-cartesien}.
Section~\ref{restrict-reflec} explains how to restrict an Olschok
model category to a full reflective
subcategory. Theorem~\ref{cartesian-reflective} encompasses all
constructions made in \cite{cubicalhdts}, \cite{homotopyprecubical}
and \cite{csts} on reflective subcategories.
Section~\ref{restrict-coreflec} explains how to restrict an Olschok
model category to a full coreflective subcategory
(Theorem~\ref{cartesian-coreflective-plus-general},
Theorem~\ref{cartesian-coreflective1} and
Theorem~\ref{cartesian-coreflective1bis}). Theorem~\ref{cartesian-coreflective-plus-general}
is implicitly used in \cite{cubicalhdts} and \cite{csts}: we prove
the statement in full generality.  Section~\ref{restrict-comma}
explains how to obtain Olschok model categories on comma categories
(Theorem~\ref{cartesian-comma}).  Finally, Section~\ref{app} is
devoted to an application of Theorem~\ref{cartesian-comma} and
Theorem~\ref{cartesian-coreflective1} to star-shaped weak transition
systems. The last section is the only one which is specific to the
theory of higher transition systems.

\subsection*{Prerequisites and notations}

All categories are locally small. The set of maps in a category $\K$
from $X$ to $Y$ is denoted by $\K(X,Y)$. The class of maps of a
category $\K$ is denoted by $\Mor(\K)$.  The composite of two maps is
denoted by $fg$ instead of $f \circ g$. The initial (final resp.)
object, if it exists, is always denoted by $\varnothing$ ($\mathbf{1}$
resp.). The identity of an object $X$ is denoted by $\id_X$.  A
subcategory is always isomorphism-closed.  Let $f$ and $g$ be two maps
of a locally presentable category $\K$. Write $f\square g$ when $f$
satisfies the \emph{left lifting property} (LLP) with respect to $g$,
or equivalently $g$ satisfies the \emph{right lifting property} (RLP)
with respect to $f$. Let us introduce the notations $\inj_\K(\C) = \{g
\in \K, \forall f \in \C, f\square g\}$ and $\cof_\K(\C) = \{f \in \K,
\forall g\in \inj_\K(\C), f\square g\}$ where $\C$ is a class of maps
of $\K$. We refer to \cite{MR95j:18001} for locally presentable
categories, and to \cite{MR2506258} for combinatorial model
categories.  We refer to \cite{MR99h:55031} and to \cite{ref_model2}
for model categories. For general facts about weak factorization
systems, see also \cite{ideeloc}. The reading of the first part of
\cite{MOPHD}, published in \cite{MO}, is recommended for any reference
about good, cartesian, and very good cylinders.

\section{Olschok model category}
\label{good-model-cat}

This is a section recalling Olschok's construction and introducing
thereby the notion of Olschok model category. Note that
Proposition~\ref{model-cartesien} is new.

\begin{nota} For every map $f:X \to Y$ and every natural transformation $\alpha : F
\to F'$ between two endofunctors of a locally presentable category
$\K$, the map $f\star \alpha$ is defined by the diagram:
\[
\xymatrix{
FX \fD{\alpha_X}\fR{f} && FY \fD{}\ar@/^15pt/@{->}[dddr]_-{\alpha_Y} &\\
&& &&\\
F'X \ar@/_15pt/@{->}[rrrd]^-{F'f}\fR{} &&  \bullet \cocartesien\ar@{->}[rd]^-{f\star \alpha} & \\
&& & F'Y.
}
\]
For a set of morphisms $\mathcal{A}$, let $\mathcal{A} \star \alpha =
\{f\star \alpha, f\in \mathcal{A}\}$.
\end{nota}

\bd \label{cyl} Let $\K$ be a locally presentable category. A {\rm
  cylinder} is a triple $(\cyl:\K \to \K,\gamma:\id \oplus
\id\Rightarrow \cyl,\sigma:\cyl \Rightarrow \id)$ consisting of a
functor $\cyl:\K \to \K$ and two natural transformations
$\gamma=\gamma^0\oplus \gamma^1 :\id \oplus \id\Rightarrow \cyl$ and
$\sigma:\cyl \Rightarrow \id$ such that the composite $\sigma 
\gamma$ is the codiagonal functor $\id \oplus \id \Rightarrow \id$ \ed

\bd \label{cyl-bon-tresbon-cartesien} Let $\K$ be a locally
presentable category. Let $(\C,\W,\F)$ be a cofibrantly generated
model structure on $\K$ where $\C$ is the class of cofibrations, $\W$
the class of weak equivalences and $\F$ the class of fibrations. A
{\rm cylinder} for $(\C,\W,\F)$ is a cylinder $(\cyl:\K \to
\K,\gamma:\id \oplus \id\Rightarrow \cyl,\sigma:\cyl \Rightarrow \id)$
such that the functorial map $\sigma_X:\cyl(X) \to X$ belongs to $\W$
for every object $X$.  The cylinder is {\rm good} if the functorial
map $\gamma_X : X \sqcup X \to \cyl(X)$ is a cofibration for every
object $X$. It is {\rm very good} if moreover the map
$\sigma_X:\cyl(X) \to X$ is a trivial fibration for every object
$X$. A good cylinder is {\rm cartesian} if
\begin{itemize}
\item The functor $\cyl:\K\to \K$ has a right adjoint $\cocyl:\K \to
  \K$ called the {\rm path functor}.
\item If $f$ is a cofibration, then so are $f\star \gamma^0$, $f\star
  \gamma^1$ and $f\star \gamma$.
\end{itemize}
\ed

The notions of Definition~\ref{cyl-bon-tresbon-cartesien} can be
adapted to a cofibrantly generated weak factorization system
$(\mathcal{L},\mathcal{R})$ by considering the combinatorial model
structure $(\mathcal{L},\Mor(\K),\mathcal{R})$. A cylinder with
respect to a set of maps $I$ is a cylinder for the weak factorization
system $(\cof_\K(I),\inj_\K(I))$, i.e. for the model structure
$(\cof_\K(I),\Mor(\K),\inj_\K(I))$.

\begin{nota} Let $I$ and $S$ be two sets of maps of a locally
  presentable category $\K$. Let $\cyl:\K\to \K$ be a cylinder with
  respect to $I$. Denote by $\Lambda_\K(\cyl,S,I)$ the set of maps
  defined as follows:
\begin{itemize}
\item $\Lambda^0_\K(\cyl,S,I) = S \cup (I \star \gamma^0) \cup (I
  \star \gamma^1)$
\item $\Lambda^{n+1}_\K(\cyl,S,I) =
  \Lambda^n_\K(\cyl,S,I) \star \gamma$
\item $\Lambda_\K(\cyl,S,I) = \bigcup_{n\geq 0}
  \Lambda^n_\K(\cyl,S,I)$.
\end{itemize} 
\end{nota}

Let us denote by $\W_\K(\cyl,S,I)$ the class of maps $f : X \to Y$ of
$\K$ such that for every object $T$ which is
$\Lambda_\K(\cyl,S,I)$-injective, the induced set map
$\K(Y,T)/\!\!\simeq \stackrel{\iso}\longrightarrow \K(X,T)/\!\!\simeq$
is a bijection, where $\simeq$ means the homotopy relation associated
with the cylinder $\cyl(-)$, i.e. for all maps $f,g:X\to Y$, $f\simeq
g$ is equivalent to the existence of a \emph{homotopy} $H:\cyl(X) \to
Y$ with $H \gamma^0 = f$ and $H \gamma^1 = g$.

\bth \label{build-model-cat} (Olschok) Let $\K$ be a locally
presentable category.  Let $I$ be a set of maps of $\mathcal{K}$. Let
$S\subset \cof_\K(I)$ be a set of maps of $\K$. Let $\cyl$ be a
cartesian cylinder for the weak factorization system
$(\cof_\K(I),\inj_\K(I))$.  Suppose that the weak factorization system
$(\cof_\K(I),\inj_\K(I))$ is cofibrant, i.e. for any object $X$ of
$\K$, the canonical map $\varnothing \to X$ belongs to
$\cof_\K(I)$. Then there exists a unique combinatorial model category
structure with class of cofibrations $\cof_\K(I)$ such that the
fibrant objects are the $\Lambda_\K(\cyl,S,I)$-injective objects. The
class of weak equivalences is $\W_\K(\cyl,S,I)$. All objects are
cofibrant.  \eth

\bpf The explanation is already given in
\cite[Theorem~2.6]{homotopyprecubical}. This theorem is a slight
modification of Olschok's main theorem \cite[Theorem~3.16]{MO} using
the characterization of fibrant objects \cite[Lemma~3.30(c)]{MO} and
the fact that a model structure is characterized by its class of
cofibrations and its class of fibrant objects:
\cite[Theorem~7.8.6]{ref_model2} works here since all objects are
cofibrant; more generally \cite[Proposition E.1.10]{quasicat} can be
used. \epf

If the cylinder is very good in Theorem~\ref{build-model-cat}, then
$\W_\K(\cyl,S,I)$ is the Grothen\-dieck localizer generated by $S$ (with
respect to the class of cofibrations $\cof_\K(I)$) by
\cite[Corollary~4.6]{MO}. In this case, $\K$ is left determined in the
sense of \cite{rotho} when $S=\varnothing$.  And the model category we
obtain for $S\neq\varnothing$ is the Bousfield localization $\bl_S\K$
of the left determined one by the set of maps
$S$.

\bp \label{model-cartesien} Let $\K$ be a combinatorial model category
such that all objects are cofibrant. Let $I$ be the set of generating
cofibrations. Let $\cyl:\K\to \K$ be a cartesian cylinder for the weak
factorization system $(\cof_\K(I),\inj_\K(I))$. Let $S\subset
\cof_\K(I)$ be a set of maps of $\K$. Then the following conditions
are equivalent:
\begin{itemize}
\item An object of $\K$ is fibrant if and only if it is
  $\Lambda_\K(\cyl,S,I)$-injective.
\item A map of $\K$ is a weak equivalence if and only if it belongs to
  $\W_\K(\cyl,S,I)$.
\end{itemize}
\ep 

\bpf Let us suppose that the fibrant objects of $\K$ are the
$\Lambda_\K(\cyl,S,I)$-injective ones. Then the model structure of
$\K$ and the one given by Theorem~\ref{build-model-cat} have the same
class of cofibrations and the same class of fibrant objects. Since all
objects are cofibrant, the class of weak equivalences is necessarily
$\W_\K(\cyl,S,I)$ by \cite[Theorem~7.8.6]{ref_model2}. Conversely, let
us suppose that a map of $\K$ is a weak equivalence if and only if it
belongs to $\W_\K(\cyl,S,I)$. Then the model structure of $\K$ and the
one given by Theorem~\ref{build-model-cat} have the same class of
cofibrations and the same class of weak equivalences. The class of
fibrations is determined by the class of trivial
cofibrations. Therefore the two model structures are equal.  So they
have the same class of fibrant objects.  \epf

\bd \label{def-good} An {\rm Olschok model category} is a
combinatorial model category satisfying the conditions of
Proposition~\ref{model-cartesien} for some cartesian cylinder $\cyl$
and some set of cofibrations $S$ called the {\rm generating anodyne
  cofibrations}.  \ed

The terminology ``anodyne'' comes from \cite{MR1924082} where the
elements of the class \[\cof_\K(\Lambda_\K(\cyl,S,I))\] are called, in
French, ``extensions anodines''. When the class of generating anodyne
cofibrations is not specified, it is supposed to be empty.

\section{Restriction to a reflective subcategory} \label{restrict-reflec}

The following theorem gives a sufficient condition for the restriction
to a full reflective subcategory of an Olschok model category to be an
Olschok model category. It implies
\cite[Theorem~9.3]{homotopyprecubical} and \cite[Theorem~5.5]{csts}
because in the latter cases the map $\eta_{\cyl(X)}$ is an
isomorphism.

\bth \label{cartesian-reflective} Let $\K$ be an Olschok model
category with generating cofibrations $I$, with generating anodyne
cofibrations $S$ and with cartesian cylinder $\cyl$. Let $\LL$ be a
full reflective locally presentable subcategory and let $\kappa:\K\to
\LL$ be the reflection. Suppose that $I=\kappa(I)$ (i.e. the source
and targets of all maps of $I$ belong to $\LL$), that $\cocyl(\LL)
\subset \LL$ where $\cocyl : \LL \to \LL$ is a right adjoint of
$\cyl:\LL \to \LL$, and that the unit map $\eta_{\cyl(X)} : \cyl(X)
\to \kappa(\cyl(X))$ has a section $s_X$ (i.e. it is split epic) for
all objects $X$ of $\LL$.  Then:
\begin{enumerate}
\item The functor $\kappa \cyl:\LL \to \LL$ is a cartesian cylinder
  with respect to $\kappa(I)$. Moreover if $\cyl:\K\to \K$ is very
  good, then $\kappa \cyl:\LL \to \LL$ is very good as well.
\item There exists a unique Olschok model structure on $\LL$ with set
  of generating cofibrations $\kappa(I)=I$, with set of generating
  anodyne cofibrations $\kappa(S)$, such that an object of $\LL$ is
  fibrant in $\LL$ if and only if it is fibrant in $\K$. The cartesian
  cylinder of $\LL$ is the functor $\kappa \cyl:\LL \to \LL$. The
  reflection $\kappa:\K \to \LL$ is a homotopically surjective (in the
  sense of \cite[Definition~3.1]{MR1870516}) left Quillen adjoint.
\end{enumerate}
\eth

Note that the existence of the section is only used to prove the
left-determinedness of the model structure of $\LL$.

\bpf By \cite[Lemma~5.2(c)]{MO}, the functor $\kappa \cyl:\LL \to \LL$
is a cartesian cylinder with respect to $\kappa(I)=I$. An object of
$\LL$ is fibrant in $\LL$ if and only if it is fibrant in $\K$ by
\cite[Lemma~5.2(b)]{MO}. The existence of the Olschok model structure
is then a consequence of Theorem~\ref{build-model-cat}.  The proof of
the fact that the reflection $\kappa:\K \to \LL$ is a homotopically
surjective left Quillen functor is mutatis mutandis the argument used
for the same fact in \cite[Theorem~9.3]{homotopyprecubical}.

\begin{figure}
\[
\xymatrix
{
A \fR{\phi} \fD{f} && \kappa(\cyl(X)) \ar@{->}[dd]^-/10pt/{\kappa(\sigma_X)} \fr{s_X}& \cyl(X) \ar@{->}[dd]^-/10pt/{\sigma_X}\fR{\eta_{\cyl(X)}} && \kappa(\cyl(X)) \ar@{->}[dd]^-/10pt/{\kappa(\sigma_X)}\\
&& & && \\
B \fR{\psi} \ar@{-->}[rruu]^-{\ell} \ar@{-->}[rrruu]^-{\ell'}&& \kappa(X)\ar@{->}[r]^-{\eta_X^{-1}} & X\ar@{->}[rr]^-{\eta_X} && \kappa(X),
}
\]
\caption{\label{zzz}$\kappa \cyl$ is very good.}
\end{figure}

Suppose now that $\cyl$ is a very good cylinder with respect to
$I$. Consider the diagram of solid arrows of $\LL$ of Figure~\ref{zzz}
where $X$ is an object of $\LL$ (this implies that $\eta_X$ is
invertible), where $f:A\to B$ belongs to $I$, and where the left-hand
square is supposed to be commutative, i.e. $\kappa(\sigma_X) \phi =
\psi f$.  The right-hand square is commutative by naturality of the
unit map of the adjunction. One has
\begin{align*}
\sigma_X  s_X &= \eta_X^{-1}  \eta_X  \sigma_X  s_X & \hbox{ since $\eta_X$ is invertible}\\
&= \eta_X^{-1}  \kappa(\sigma_X)  \eta_{\cyl(X)}  s_X & \hbox{ by naturality of the unit map}\\
&= \eta_X^{-1}  \kappa(\sigma_X) & \hbox{ by hypothesis on $s_X$.}
\end{align*}
This means that the middle square is commutative as well.  One deduces
that the composite of the left-hand square and the middle square is a
commutative square, i.e. $\sigma_X s_X \phi = \eta_X^{-1} \psi
f$. Since $\cyl$ is a very good cylinder of $\K$ with respect to $I$,
there exists a lift $\ell':B\to \cyl(X)$ such that $\sigma_X \ell' =
\eta_X^{-1} \psi$ and $\ell' f = s_X \phi$.  Let $\ell =
\eta_{\cyl(X)} \ell'$. One has
\begin{align*}
\kappa(\sigma_X)  \ell &= \kappa(\sigma_X)  \eta_{\cyl(X)}  \ell' & \hbox{ by definition of $\ell$}\\
&= \eta_X  \sigma_X  \ell' & \hbox{ since the right-hand square is commutative}\\
&=\eta_X  \eta_X^{-1}  \psi & \hbox{ by definition of $\ell'$}\\
&=\psi & \hbox{ by trivial simplification}.
\end{align*}
And one has 
\begin{align*}
\ell f &= \eta_{\cyl(X)}  \ell'  f & \hbox{ by definition of $\ell$}\\
&= \eta_{\cyl(X)}  s_X  \phi & \hbox{ by definition of $\ell'$} \\
&= \phi &\hbox{ by hypothesis on $s_X$.}
\end{align*}
Therefore $\ell$ is a lift for the left-hand square. Hence the
cylinder $\kappa \cyl:\LL\to \LL$ is very good with respect to
$I$. The proof is complete.  \epf

\begin{cor}
  With the notations of Theorem~\ref{cartesian-reflective}, there
  exists a Bousfield localization of $\K$ which is Quillen equivalent
  to $\LL$.
\end{cor}

\begin{cor} With the notations of Theorem~\ref{cartesian-reflective},
  the inclusion $\LL\subset \K$ reflects weak equivalences.
\end{cor}

\bpf Let $f:X\to Y$ be a map of $\LL$ which is a weak equivalence of
$\K$.  Then for any fibrant object $F$ of $\K$, the set map
$\K(Y,F)\to \K(X,F)$ induced by composing by $f$ gives rise to a
bijection between the homotopy classes. Since the fibrant objects of
$\LL$ are the fibrant objects of $\K$ belonging to $\LL$, this implies
that $f$ is a weak equivalence of $\LL$.  \epf

\section{Restriction to a coreflective subcategory} \label{restrict-coreflec}

The following theorem is the general theorem behind the construction
of the homotopy theory of cubical transition systems in
\cite{cubicalhdts}.

\bth \label{cartesian-coreflective-plus-general} Let $\K$ be an
Olschok model category with cartesian cylinder $\cyl$ and with set of
generating cofibrations $I$. Let $\LL$ be a full coreflective locally
presentable subcategory such that:
\begin{itemize}
\item There exists a set of maps $J$ such that $\cof_\LL(J) =
  \cof_\K(I)\cap \Mor(\LL)$.
\item $\cyl(\LL)\subset \LL$.
\end{itemize}
Then there exists a structure of Olschok model category on $\LL$ such
that the cofibrations are the cofibrations of $\K$ between objects of
$\LL$ and such that the restriction to $\LL$ of $\cyl$ is a cartesian
cylinder for this model structure. Moreover, if $\cyl$ is very good in
$\K$, then its restriction to $\LL$ gives rise to a very good cylinder 
in $\LL$.  \eth

\bpf The set $J$ will be the set of generating cofibrations of the
Olschok model category $\LL$. Let $A$ be an object of $\LL$. Consider
the factorization of the codiagonal of $A$ given by this cylinder:
\[\xymatrix{A\sqcup A \fR{\gamma_A} && \cyl(A) \fR{\sigma_A} && A.}\]
By hypothesis, $\cyl(A)$ is an object of $\LL$. Therefore $\gamma_A$
is a cofibration of $\LL$. So the restriction of $\cyl$ to $\LL$ gives
rise to a good cylinder. Let $f:A\to B$ be a cofibration of
$\LL$. Then the maps $f\star \gamma^\epsilon : B \sqcup_A \cyl(A)
\longrightarrow B$ for $\epsilon=0,1$ and $f\star \gamma : (B\sqcup B)
\sqcup_{A\sqcup A} \cyl(A) \longrightarrow B \sqcup B$ are
cofibrations of $\K$ since $\cyl$ is a cartesian cylinder. The sources
and the targets of these maps belong to $\LL$ since $\LL$ is a
coreflective subcategory. So the maps $f\star \gamma^\epsilon$ for
$\epsilon=0,1$ and $f\star \gamma$ are cofibrations of $\LL$. Let $A$
and $B$ be two objects of $\LL$. Then 
\begin{align*}
\LL(\cyl(A),B) & = \K(\cyl(A),B) & \hbox{since $\LL$ is a full
  subcategory} \\ & = \K(A,\cocyl(B)) & \hbox{where $\cocyl$ is a
  right adjoint of $\cyl$} \\ & = \LL(A,\xi(\cocyl(B))) & \hbox{where
  $\xi$ is the coreflection.}
\end{align*} 
This implies that the restriction of $\cyl$ to $\LL$ gives rise to a
cartesian cylinder. The proof of the existence of the model structure
is complete thanks to Theorem~\ref{build-model-cat}.

Let us suppose now that $\cyl$ is very good in $\K$. Then for every
object $A$ of $\LL$, the map $\sigma_A : \cyl(A)\to A$ is a trivial
fibration of $\K$ which satisfies the RLP with respect to any
cofibration of $\K$. Since the cofibrations of $\LL$ are exactly the
cofibrations of $\K$ between objects of $\LL$, the map $\sigma_A :
\cyl(A)\to A$ is a trivial fibration of $\LL$ as well.  \epf

\bth \label{cartesian-coreflective2} With the notations and hypotheses
of Theorem~\ref{cartesian-coreflective-plus-general}, assume that the
set $S$ of generating anodyne cofibrations of $\K$ belongs to
$\LL$. Let us equip $\LL$ with the Olschok model structure having the
same set of generating anodyne cofibrations. Then the inclusion
functor $\LL \to \K$ is a left Quillen functor. \eth

\bpf It is mutatis mutandis the proof of
\cite[Theorem~6.3]{cubicalhdts}.  \epf

Theorem~\ref{cartesian-coreflective-plus-general} has the following
corollaries:

\bth \label{cartesian-coreflective1} Let $\K$ be an Olschok model
category with cartesian cylinder $\cyl$. Let $\LL$ be a full
coreflective subcategory such that:
\begin{itemize}
\item $\LL$ is a small cone-injectivity class with respect to a set of
  cofibrations of $\K$.
\item $\cyl(\LL)\subset \LL$.
\end{itemize}
Then there exists a structure of Olschok model category on $\LL$ such
that the cofibrations are the cofibrations of $\K$ between objects of
$\LL$ and such that the restriction to $\LL$ of $\cyl$ is a cartesian
cylinder for this model structure. Moreover, if $\cyl$ is very good in
$\K$, then its restriction to $\LL$ gives rise to a very good cylinder
in $\LL$.  \eth

Note that this is the theorem used in \cite{cubicalhdts}.

\bpf Since $\LL$ is full coreflective, it is cocomplete. And since it
is a small cone-injectivity class, it is accessible by
\cite[Proposition~4.16]{MR95j:18001}. Therefore $\LL$ is locally
presentable.  Let $I$ be the set of generating cofibrations of
$\K$. By \cite[Theorem~A.5]{cubicalhdts}, there exists a set of maps
$J$ such that $\cof_\LL(J) = \cof_\K(I)\cap \Mor(\LL)$. We can then
apply Theorem~\ref{cartesian-coreflective-plus-general}. \epf

\bth \label{cartesian-coreflective1bis} Let $\K$ be an Olschok model
category with cartesian cylinder $\cyl$ with set of generating
cofibrations $I$. Let $\LL$ be a full coreflective locally presentable
subcategory such that:
\begin{itemize}
\item $I$ has a solution set $J \subset \cof_\K(I)$ with respect to
  $\LL$, i.e. $J$ is a set of maps of $\LL$ such that every map $i
  \rightarrow w$ of $\Mor(\K)$ from $i \in I$ to $w\in \Mor(\LL)$
  factors as a composite $i \rightarrow j \rightarrow w$ with $j\in
  J$.
\item $\cyl(\LL)\subset \LL$.
\end{itemize}
Then there exists a structure of Olschok model category on $\LL$ such
that the cofibrations are the cofibrations of $\K$ between objects of
$\LL$ and such that the restriction to $\LL$ of $\cyl$ is a cartesian
cylinder for this model structure. Moreover, if $\cyl$ is very good in
$\K$, then its restriction to $\LL$ gives rise to a very good cylinder
in $\LL$.  \eth

\bpf By \cite[Lemma~A.3]{cubicalhdts}, there is the equality
$\cof_\LL(J) = \cof_\K(I)\cap \Mor(\LL)$. We can then apply
Theorem~\ref{cartesian-coreflective-plus-general}.  \epf

\section{Olschok model category and comma category}
\label{restrict-comma}

The following well-known proposition introduces some useful notations:

\bp Let $\K$ be a locally presentable category. Let $i$ be an object
of $\K$.  The forgetful functor $\omega^i:i\ddownarrow \K \to \K$
defined on objects by $\omega^i(i\to X) = X$ and on maps by
$\omega^i(i\to f)=f$ is a right adjoint. In particular, it is
limit-preserving. A colimit in the comma category $i\ddownarrow \K$ is
obtained by taking the colimit in $\K$ of the cone with top the object
$i$ and with basis the diagram of underlying objects of $\K$.  The
forgetful functor $\omega^i:i\ddownarrow \K \rightarrow \K$ commutes
with colimits of connected diagrams (and in particular, it is
accessible).  \ep

Note that the forgetful functor $\omega^i:(i\ddownarrow \K)
\rightarrow \K$ does not preserve binary coproducts. Indeed, the
binary coproduct of $i\to X$ and $i\to Y$ is the amalgamated sum $i\to
X\sqcup_{i} Y$. 

\bpf The left adjoint $\rho^i:\K \rightarrow i\ddownarrow \K$ is
defined on objects by $\rho^i(X) = (i\to i \sqcup X)$ and on morphisms
by $\rho^i(f) = \id_{i} \sqcup f$. The last assertions are clear. \epf

Let $\K$ be a locally presentable category. Let $i$ be an object of
$\K$. Then the comma category $i\ddownarrow \K$ is locally presentable
by \cite[Proposition~1.57]{MR95j:18001}. Let $(\C,\W,\F)$ be a
cofibrantly generated model structure on $\K$.  Then the triple
\[((\omega^i)^{-1}(\C), (\omega^i)^{-1}(\W), (\omega^i)^{-1}(\F))\] is a
cofibrantly generated model structure on $i\ddownarrow \K$ by
\cite[Theorem~2.7]{undercat}. If $I$ is the set of generating
cofibrations of $\K$, then the set of generating cofibrations of the
comma category $i\ddownarrow \K$ is the set $\rho^i(I)$ where
$\rho^i:\K \to i\ddownarrow \K$ is the left adjoint of the functor
$\omega^i$ above defined.

\begin{lem} \label{l0} Let $\cyl:\K\to \K$ be a cylinder functor of a
  locally presentable category $\K$. Assume that it has a right
  adjoint $\cocyl:\K\to\K$. Let $i$ be an object of $\K$. Define the
  functor $\cyl_i:i\ddownarrow \K \to i\ddownarrow \K$ by the natural
  pushout diagram
\[
\xymatrix
{
\cyl(i) \fR{\sigma_i}\fD{} && i  \ar@{->}[dd]^-{\cyl_i(i\to X)} \\
&& \\
\cyl(X)  \fR{} && \omega^i(\cyl_i(i\to X)) \cocartesien 
}
\]
for every object $X$ of $\K$ and $\cocyl_i:i\ddownarrow \K \to
i\ddownarrow \K$ by the natural diagram \[\cocyl_i(i\to Y) := i\longrightarrow \cocyl(i)
\longrightarrow \cocyl(Y)\] for every object $i\to
Y$ of $i\ddownarrow \K$ where $i\longrightarrow \cocyl(i)$ is the map
corresponding to $\sigma_i:\cyl(i)\to i$ by the adjunction.  Then
$\cyl_i:i\ddownarrow \K \to i\ddownarrow \K$ is left adjoint to
$\cocyl_i:i\ddownarrow \K \to i\ddownarrow \K$.
\end{lem}

Note that it can be easily checked that the functor
$\cocyl_i:i\ddownarrow \K \to i\ddownarrow \K$ is accessible and
limit-preserving. Therefore, by \cite[Theorem~1.66]{MR95j:18001}, it
has a left adjoint since the category $i\ddownarrow \K$ is locally
presentable.

\bpf Let $i\to X$ and $i\to Y$ be fixed. By definition of $\cyl_i$,
there is a bijection between the sets of commutative diagrams
\[
\left\{
\begin{CD}
i @=   i\\
@VVV @VVV \\
\omega^i(\cyl_i(i\to X)) @>>> Y
\end{CD}
\right\} \iso 
\left\{
\begin{CD}
\cyl(i)  @>>>   i\\
@VVV @VVV \\
\cyl(X) @>>> Y
\end{CD}
\right\}.
\]
By adjunction, there is a bijection between the sets of commutative
diagrams
\[
\left\{
\begin{CD}
\cyl(i)  @>>>   i\\
@VVV @VVV \\
\cyl(X) @>>> Y
\end{CD}
\right\}
\iso
\left\{
\begin{CD}
i  @>>>   \cocyl(i)\\
@VVV @VVV \\
X @>>> \cocyl(Y)
\end{CD}
\right\}.
\]
Finally, by the definition of $\cocyl_i$, there is a bijection between
the sets of commutative diagrams
\[
\left\{
\begin{CD}
i  @>>>   \cocyl(i)\\
@VVV @VVV \\
X @>>> \cocyl(Y)
\end{CD}
\right\}
\iso
\left\{
\begin{CD}
i  @=   i\\
@VVV @VVV \\
X @>>> \omega^i(\cocyl_i(Y))
\end{CD}
\right\}.
\]
\epf

\begin{lem} \label{l1} Let $\K$ be a locally presentable category. Let
  $i$ be an object of $\K$. Let $s:A\to B$ be a map of $\K$. Let $i\to
  X$ be an object of the comma category $i\ddownarrow \K$. Then $i\to
  X$ is injective with respect to $\rho^i(s): i\sqcup A \to i\sqcup B$
  if and only if $X = \omega^i(i\to X)$ is injective with respect to
  $s$. \end{lem}

\bpf One has the commutative diagram of sets:
\[
\xymatrix
{
\K(B,X) \fD{f\mapsto f s} \ar@{=}[rr] && \fD{f\mapsto f s} \K(B,\omega^i(i\to X))\fR{\iso} && (i\ddownarrow \K)(\rho^i(B),X) \fD{g\mapsto g \rho^i(s)}\\
&& && \\
\K(A,X)  \ar@{=}[rr] &&  \K(A,\omega^i(i\to X))\fR{\iso} && (i\ddownarrow \K)(\rho^i(A),X).
}
\]
Therefore the left vertical arrow is onto if and only if the right vertical arrow is onto as well. 
\epf

\begin{cor} \label{l2} Let $\Lambda$ be a set of maps of a locally
  presentable category $\K$.  Let $i$ be an object of $\K$. Then an
  object $i\to X$ of $i\ddownarrow \K$ is $\rho^i(\Lambda)$-injective
  if and only if $X$ is $\Lambda$-injective.
\end{cor}

\begin{lem} \label{l3} With the notations and hypotheses of
  Lemma~\ref{l0}, let $A$ be an object of $\K$. Then there is the
  natural isomorphism $\cyl_i(\rho^i(A)) \iso \rho^i(\cyl(A))$.
\end{lem}

\bpf One has the bijections
\begin{align*}
(i\ddownarrow \K)(\cyl_i(\rho^i(A)),i\to B) & \iso (i\ddownarrow \K)(\rho^i(A),\cocyl_i(i \to B)) & \hbox { by adjunction}\\
 & \iso \K(A,\omega^i(\cocyl_i(i \to B))) & \hbox { by adjunction}\\
& \iso \mathcal{K}(A,\cocyl(B)) & \hbox{ by definition of $\cocyl_i$}\\
& \iso \mathcal{K}(\cyl(A),B) & \hbox{ by adjunction}\\
& \iso \mathcal{K}(\cyl(A),\omega^i(i\to B)) & \hbox{ by definition of $\omega^i$}\\
& \iso (i\ddownarrow \K)(\rho^i(\cyl(A)),i\to B) & \hbox{ by adjunction again.}
\end{align*}
Hence the result by Yoneda.
\epf

\begin{lem} \label{l4} Let $\K$ be a locally presentable category. Let
  $i\to X$ be an object of $i\ddownarrow \K$. Then one has the pushout
  diagrams
\[
\xymatrix{
i\sqcup i \fR{\id_i\sqcup\id_i} \fD{} && i \fD{} &&  i\sqcup i \fR{\id_i\sqcup\id_i} \fD{} && i \fD{}\\
&& && &&\\
X \sqcup X \fR{} && X\sqcup_i X, \cocartesien&& X  \ar@{=}[rr] && X. \cocartesien
}
\]
\end{lem}

\bpf
Consider the pushout diagrams of $\K$: 
\[
\xymatrix{
i\sqcup i \fR{\id_i\sqcup\id_i} \fD{} && i \fD{} &&  i\sqcup i \fR{\id_i\sqcup\id_i} \fD{} && i \fD{}\\
&& && &&\\
X \sqcup X \fR{} && Z, \cocartesien&& X  \fR{} && T. \cocartesien
}
\]
Let $U$ be an object of $\K$. One has the pullback diagram of sets
\[
\xymatrix{
\K(Z,U) \cartesien\fR{} \fD{} && \K(i,u) \fD{} \\
&& && &&\\
\K(X,U)\p \K(X,U) \fR{} && \K(i,U)\p \K(i,U).
}
\]
Therefore one obtains the bijections of sets 
\begin{multline*}
\K(Z,U) \iso (\K(X,U)\p \K(X,U)) \p_{\K(i,U)\p \K(i,U)} \K(i,U) \\ \iso \K(X,U)\p_{\K(i,U)} \K(X,U) \iso \K(X\sqcup_i X,U).
\end{multline*}
By Yoneda, one obtains the isomorphism $Z\iso X\sqcup_i X$. And one
has the pullback of sets 
\[
\xymatrix{
 \K(T,U) \cartesien\fR{} \fD{} && \K(i,U) \fD{}\\
 &&\\
\K(X,U)  \fR{} && \K(i,U)\p \K(i,U). 
}
\]
Therefore one obtains the bijections of sets \[\K(T,U) \iso \K(X,U)
\p_{\K(i,U)\p \K(i,U)}\K(i,U) \iso \K(X,U).\] By Yoneda, one obtains
the isomorphism $T\iso X$.  \epf

\begin{lem} \label{l0bis} Let $\cyl:\K\to \K$ be a cylinder functor of
  a locally presentable category $\K$. Assume that it has a right
  adjoint $\cocyl:\K\to\K$. Let $i$ be an object of $\K$ such that the
  map $\gamma_i:i\sqcup i\to \cyl(i)$ is epic. Then there is a pushout
  diagram
\[
\xymatrix
{
i \sqcup i \fR{} \fd{} && i  \ar@{->}[dd]^-{\cyl_i(i\to X)} \\
X \sqcup X \fd{\gamma_X}  && \\
\cyl(X)  \fR{} && \omega^i(\cyl_i(i\to X)) \cocartesien 
}
\]
in $\K$ for every object $i\to X$ of $i\ddownarrow \K$.
\end{lem}

\begin{figure}
\[
\xymatrix
{
i\sqcup i \ar@{->}[rrrr]^-{\id_i\sqcup \id_i} \fD{e_X\sqcup e_X} \ar@{->}[rrdd]^-{\gamma_i} && && i \ar@{=}[dd]\\
&& && \\
X\sqcup X \fD{\gamma_X} && \cyl(i) \fD{\cyl(e_X)} \fR{\sigma_i} && i \fD{}\\
&& && \\
\cyl(X) \ar@{=}[rr] && \cyl(X) \fR{} && Y.
}
\]
\caption{\label{yyy}Isomorphism between two categories of commutative diagrams.}
\end{figure}
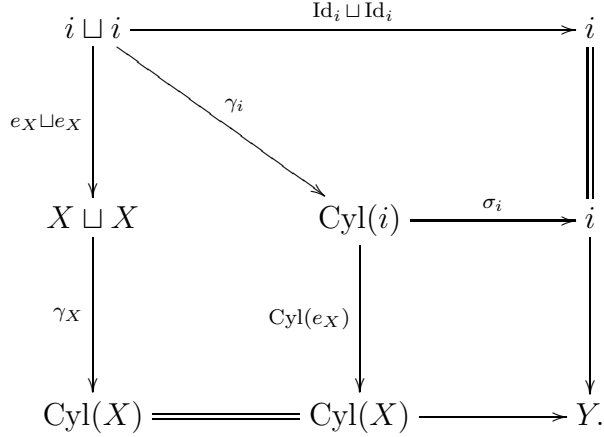

\bpf Let $e_X:i\to X$ be a fixed object of $i\ddownarrow \K$.
Consider a diagram of the form of Figure~\ref{yyy}.  We obtain a map
$F$ between the set of squares
\[
\left\{
\begin{CD}
\cyl(i)  @>>>   i\\
@VVV @VVV \\
\cyl(X) @>>> Y
\end{CD}
\right\}
\stackrel{F}\longrightarrow
\left\{
\begin{CD}
i\sqcup i  @>>>   i\\
@VVV @VVV \\
\cyl(X) @>>> Y
\end{CD}
\right\}
.
\]
If $D$ is a commutative square, then $F(D)$ is a commutative square.
Since the map $\gamma_i:i\sqcup i \to \cyl(i)$ is epic, if $F(D)$ is a
commutative square, then $D$ is a commutative square as well: $D$ is
commutative if and only if $F(D)$ is commutative. We have obtained a
bijection between the sets of commutative diagrams
\[
\left\{
\begin{CD}
\cyl(i)  @>>>   i\\
@VVV @VVV \\
\cyl(X) @>>> Y
\end{CD}
\right\}
\iso
\left\{
\begin{CD}
i\sqcup i  @>>>   i\\
@VVV @VVV \\
\cyl(X) @>>> Y
\end{CD}
\right\}
\]
which gives rise to an isomorphism between the corresponding
categories of commutative diagrams. The initial objects are the
pushout diagrams.  \epf

The main theorem of this section is the following one:

\bth \label{cartesian-comma} Let $\K$ be an Olschok model category
with the set of generating cofibrations $I$, the set of generating
anodyne cofibration $S$, and the cartesian cylinder
$\cyl:\K\to\K$. Let $i$ be an object of $\K$ such that every map of
$\K$ with source $i$ is a cofibration and such that the map
$\gamma_i:i\sqcup i \to \cyl(i)$ is epic.  Then the combinatorial
model category $i\ddownarrow \K$ is Olschok as well. The set of
generating cofibrations of $i\ddownarrow \K$ is $\rho^i(I)$. The set
of generating anodyne cofibrations of $i\ddownarrow \K$ is
$\rho^i(S)$. The cartesian cylinder is the functor
$\cyl_i:i\ddownarrow \K \to i\ddownarrow \K$ defined in
Lemma~\ref{l0}. \eth

Note that the condition ``$\gamma_i:i\sqcup i \to \cyl(i)$ epic'' is
not satisfied for the model category of topological spaces: the
inclusion map $\{0,1\} \subset [0,1]$ is not an epimorphism. We will
see an example of such a situation in Section~\ref{app}. Other
examples of such a situation can be obtained by using the category of
labelled symmetric precubical sets \cite{homotopyprecubical}, the
category of flows \cite{model3} or the category of multipointed
$d$-spaces \cite{mdtop} with $i=\{0\}$. We have for all these examples
$\cyl(i)=i$. The map $\gamma_i:i\sqcup i \to \cyl(i)$ is then the
epimorphism $R:\{0,1\} \to \{0\}$. Note that the model categories of
topological spaces, of flows and of multipointed $d$-spaces are not
Olschok model categories since they contain non-cofibrant objects. But
it can be proved that they are left determined. The model category of
labelled symmetric precubical sets of \cite{homotopyprecubical} is an
Olschok model category. However, it is not known if the latter is left
determined.

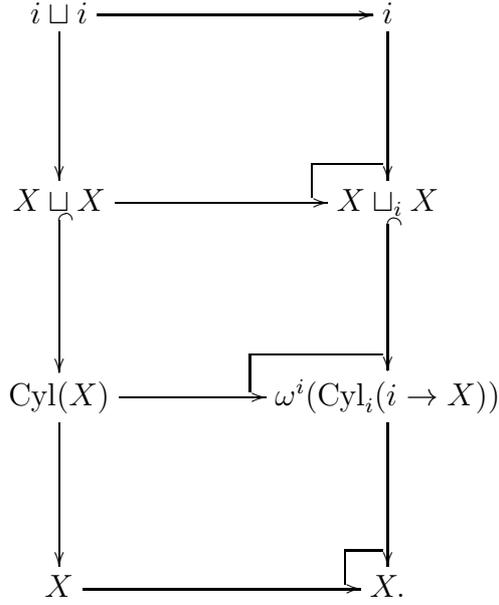
\begin{figure}
\[
\xymatrix
{
i \sqcup i \fR{} \fD{} && i\fD{} \\
&& \\
X \sqcup X \fR{} \ar@{^(->}[dd] && X \sqcup_i X  \ar@{^(->}[dd] \cocartesien  \\
&& \\
\cyl(X)  \fR{} \fD{} && \omega^i(\cyl_i(i\to X)) \cocartesien \fD{} \\
&& \\
X  \fR{} && X. \cocartesien
}
\]
\caption{\label{ttt}Composite of three pushout squares (under the hypothesis $\gamma_i:i \sqcup i \to \cyl(i)$ epic).}
\end{figure}

\begin{figure}
\[
\xymatrix
{
&& i\sqcup i \fR{} \fD{} && i \ar@{->}[dd]^-{\cyl_i(i\to X)} \\
&& && \\
X \ar@{^(->}[rr]^-{\gamma^\epsilon_X} \ar@{^(->}[dd] && \cyl(X)\ar@{^(->}[dd] \fR{} &&  \omega^i(\cyl_i(i\to X))  \ar@{^(->}[dd] \cocartesien \\
&& && \\
Y \ar@{^(->}[rr]^-{} \ar@{=}[dd] && \bullet\ar@{^(->}[dd]^-{\omega_i(f)\star \gamma^\epsilon} \fR{} \cocartesien && \bullet \ar@{^(->}[dd]^-{} \cocartesien\\
&& && \\
Y \ar@{^(->}[rr]^-{\gamma^\epsilon_Y} && \cyl(Y)\fR{} && Z, \cocartesien
}
\]
\caption{\label{xxx}$\cyl_i$ is cartesian.}
\end{figure}
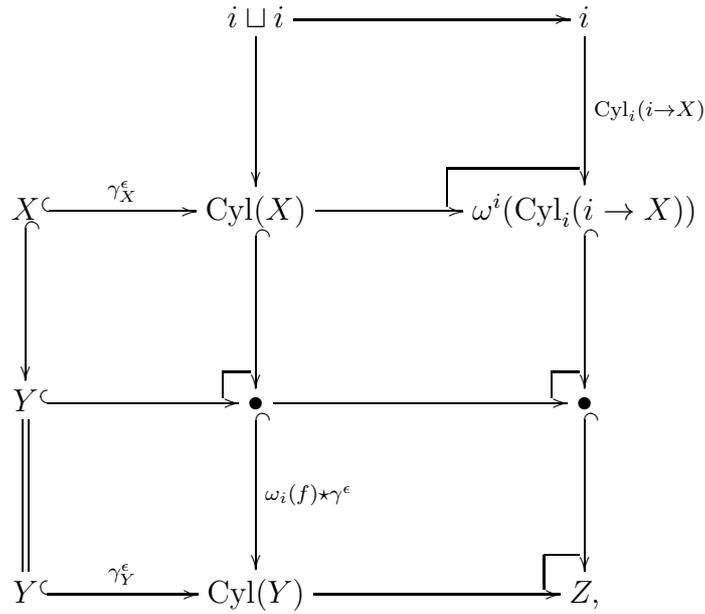

\bpf Since every map with source $i$ is a cofibration and since the
identity of $i$ is the initial object of $i\ddownarrow \K$, all
objects of the model category $i\ddownarrow \K$ are cofibrant. Let
$i\to X$ be an object of $i\ddownarrow \K$. Consider the composite
diagram of $\K$ of Figure~\ref{ttt}.  By Lemma~\ref{l0bis} and
Lemma~\ref{l4}, the three squares above are pushout squares: in
particular, they are commutative.  The commutativity of
Figure~\ref{ttt} implies that the natural map $X\sqcup_i X \to
\cyl_i(X) \to X$ is the codiagonal of $i\to X$ in $i\ddownarrow \K$.
Since the functor $\cyl:\K\to \K$ is a good cylinder, the map
$X\sqcup_i X \to \omega^i(\cyl_i(i\to X))$ is a cofibration of
$\K$. Therefore the functor $\cyl_i:i\ddownarrow \K \to i\ddownarrow
\K$ is a good cylinder for the set of maps $\rho^i(I)$~\footnote{Note
  that we do not know yet that the map $\cyl_i(X) \to X$ is a weak
  equivalence; this fact will be a consequence of this theorem. So we
  cannot yet say that $\cyl_i:i\ddownarrow \K \to i\ddownarrow \K$ is
  a good cylinder for the model category $i\ddownarrow \K$.}. By
Lemma~\ref{l0}, the functor $\cyl_i:i\ddownarrow \K \to i\ddownarrow
\K$ has a right adjoint.  Let $f:X\to Y$ be a cofibration of the comma
category $i\ddownarrow \K$. Then $\omega^i(f)$ is a cofibration by
definition of the model category $i\ddownarrow \K$. One has the
commutative diagram of $\K$ of Figure~\ref{xxx}, with $\epsilon=0,1$
where $Z$ is defined as the pushout of the right-bottom square.  Since
we have the pushout diagram
\[
\xymatrix
{
i\sqcup i \fR{} \fD{} && i \fD{} \\
&& \\
\cyl(Y) \fR{} && \cocartesien Z,
}
\]
one deduces that $\omega^i(\cyl_i(i\to Y)) =Z$ and we obtain the pushout diagram
\[
\xymatrix
{
\bullet \ar@{^(->}[dd]^-{\omega_i(f)\star \gamma^\epsilon} \fR{} && \bullet  \ar@{^(->}[dd]^-{f\star \gamma^\epsilon} \\
&& \\
\cyl(Y) \fR{} && \omega^i(\cyl_i(i\to Y)).
}
\]
Therefore the map $f\star \gamma^\epsilon$ is a cofibration of the
comma category $i\ddownarrow \K$.  We prove in the same way that
$f\star\gamma$ is a cofibration. Hence the functor
$\cyl_i:i\ddownarrow \K \to i\ddownarrow \K$ is a cartesian cylinder
for $\rho^i(I)$.  By Theorem~\ref{build-model-cat}, we deduce that
there exists a unique Olschok model category structure on
$i\ddownarrow \K$ with the set of generating cofibrations $\rho^i(I)$,
with the set of generating anodyne cofibrations $\rho^i(S)$, with the
cartesian cylinder $\cyl_i:i\ddownarrow \K \to i\ddownarrow \K$ and
such that an object is fibrant if and only if it is
$\Lambda_{i\ddownarrow \K}(\cyl_i,\rho^i(S),\rho^i(I))$-injective.

Let $f:A\to B$ be a map of $\K$. Since the functor $\rho^i:\K \to
i\ddownarrow \K$ preserves colimits, one has the commutative diagram
of solid arrows of $i\ddownarrow \K$
\[
\xymatrix{
\rho^i(A) \fD{\rho^i(\gamma^\epsilon_A)}\fR{\rho^i(f)} && \rho^i(B) \fD{\rho^i(\gamma^\epsilon_B)}\ar@/^20pt/@{->}[dddr] &\\
&& &&\\
\rho^i(\cyl(A)) \ar@/_20pt/@{->}[rrrd]_-{\rho^i(\cyl(f))}\fR{} &&  \bullet \cocartesien\ar@{->}[rd]^-{\rho^i(f\star \gamma^\epsilon)} & \\
&& & \rho^i(\cyl(B))
}
\]
for $\epsilon=0,1$. By Lemma~\ref{l3}, one deduces that $\rho^i(f)
\star \gamma^\epsilon = \rho^i(f \star \gamma^\epsilon)$ for $\epsilon
= 0,1$. For the same reason, one has the commutative diagram of solid
arrows of $i\ddownarrow \K$
\[
\xymatrix{
\rho^i(A)\sqcup \rho^i(A) \fD{\rho^i(\gamma_A)}\fR{\rho^i(f)} && \rho^i(B)\sqcup \rho^i(B) \fD{\rho^i(\gamma_B)}\ar@/^20pt/@{->}[dddr] &\\
&& &&\\
\rho^i(\cyl(A)) \ar@/_20pt/@{->}[rrrd]_-{\rho^i(\cyl(f))}\fR{} &&  \bullet \cocartesien\ar@{->}[rd]^-{\rho^i(f\star \gamma)} & \\
&& & \rho^i(\cyl(B)).
}
\]
By Lemma~\ref{l3}, one deduces that $\rho^i(f) \star \gamma = \rho^i(f
\star \gamma)$. So by Corollary~\ref{l2}, an object $i\to X$ of the
comma category $i\ddownarrow \K$ is $\Lambda_{i\ddownarrow
  \K}(\cyl_i,\rho^i(S),\rho^i(I))$-injective if and only if $X$ is
$\Lambda_{\K}(\cyl,S,I)$-injective, i.e. if and only if $X$ is fibrant
in $\K$.

We deduce that the model category constructed in this proof has the
same cofibrations and the same fibrant objects as the model category
$i\ddownarrow \K$. Hence they are equal by
\cite[Theorem~7.8.6]{ref_model2} since all objects are
cofibrant~\footnote{Moreover, we can say now that the map $\cyl_i(X)
  \to X$ is a weak equivalence of $\K$ as well, which was not possible
  earlier.}. \epf

It is not clear how to prove without additional hypothesis that if
$\cyl$ is very good, then $\cyl_i$ is very good as well. In the
situations one wants to use this construction, the map $\cyl(X) \to
\cyl_i(X)$ is always split epic.  In this case, one has:

\begin{cor}\label{far-fetched}  With the same notations and hypotheses as in
  Theorem~\ref{cartesian-comma}, if the map $p_X:\cyl(X) \to
  \omega^i(\cyl_i(X))$ is split epic for every $X$, then if $\cyl$ is
  a very good cylinder of $\K$, then $\cyl_i$ is a very good cylinder
  of $i\ddownarrow \K$.
\end{cor}

\bpf We start from a commutative diagram of $\K$ where $f$ is a map of $I$: 
\[
\xymatrix
{
&& \cyl(X) \fD{p_X}\ar@/^50pt/[dddd]^-{\sigma_X}\\
&& \\
i \sqcup A \ar@{-->}[rruu]^-{k}\fR{\phi}\fD{\rho^i(f)} &&  \omega^i(\cyl_i(X)) \fD{\sigma^i_X} \\
&& \\
i \sqcup B \fR{\psi} && X.
}
\]
Let $k=s_X \phi$ where $s_X:\omega^i(\cyl_i(X)) \to \cyl(X)$ is a
section of the split epic $\cyl(X)\to \omega^i(\cyl_i(X))$. Since $i$
is cofibrant, $\rho^i(f)$ is a cofibration of $\K$. Since $\cyl$ is
very good, there exists a lift $\ell:i\sqcup B \to \cyl(X)$ such that
$\ell \rho^i(f) = k$ and $\sigma_X \ell = \psi$. Then one has $(p_X
\ell) \rho^i(f) = p_X k = p_X s_X \phi = \phi$ and $\sigma^i_X (p_X
\ell) = \sigma_X \ell= \psi$. Hence the cylinder $\cyl_i$ is a very
good cylinder of $i\ddownarrow \K$.  \epf

\section{The homotopy theory of star-shaped weak transition systems}
\label{app}

Weak transition systems are introduced in \cite{hdts} as a rewording
of Cattani-Sassone's notion of higher dimensional transition system
\cite{MR1461821}. The purpose of these combinatorial objects is to
model the concurrent execution of $n$ actions by a transition between
two states labelled by a multiset $\{u_1,\dots,u_n\}$ of actions. The
category of weak transition systems is a convenient category to study
these objects from a categorical and homotopical point of view
\cite{hdts} \cite{cubicalhdts} \cite{csts}.

\begin{nota} Let $\Sigma$ be a fixed nonempty set of {\rm labels}.  \end{nota}

\bd\label{def_HDTS} A {\rm weak transition system} consists of a
triple \[X=(S,\mu:L\rightarrow \Sigma,T=\bigcup_{n\geq 1}T_n)\] where
$S$ is a set of {\rm states}, where $L$ is a set of {\rm actions},
where $\mu:L\rightarrow \Sigma$ is a set map called the {\rm labelling
  map}, and finally where $T_n\subset S\p L^n\p S$ for $n \geq 1$ is a
set of {\rm $n$-transitions} or {\rm $n$-dimensional transitions} such
that one has:
\begin{itemize}
\item (Multiset axiom) For every permutation $\sigma$ of
  $\{1,\dots,n\}$ with $n\geq 2$, if the tuple
  $(\alpha,u_1,\dots,u_n,\beta)$ is a transition, then the tuple
  $(\alpha,u_{\sigma(1)}, \dots, u_{\sigma(n)}, \beta)$ is a
  transition as well.
\item (Patching axiom~\footnote{This axiom is called the Coherence
    axiom in \cite{hdts} and \cite{cubicalhdts}, and the composition
    axiom in \cite{csts}.}) For every $(n+2)$-tuple
  $(\alpha,u_1,\dots,u_n,\beta)$ with $n\geq 3$, for every $p,q\geq 1$
  with $p+q<n$, if the five tuples 
\begin{align*}
 &(\alpha,u_1, \dots, u_n, \beta), \\ &(\alpha,u_1, \dots, u_p,
  \nu_1), (\nu_1, u_{p+1}, \dots, u_n, \beta),\\ &(\alpha, u_1,
  \dots, u_{p+q}, \nu_2), (\nu_2, u_{p+q+1}, \dots, u_n, \beta)
\end{align*}
  are transitions, then the $(q+2)$-tuple $(\nu_1, u_{p+1}, \dots,
  u_{p+q}, \nu_2)$ is a transition as well.
\end{itemize}
A map of weak transition systems
\[f:(S,\mu : L \rightarrow \Sigma,(T_n)_{n\geq 1}) \rightarrow
(S',\mu' : L' \rightarrow \Sigma ,(T'_n)_{n\geq 1})\] consists of a
set map $f_0: S \rightarrow S'$ and a commutative square
\[
\xymatrix{
  L \ar@{->}[r]^-{\mu} \ar@{->}[d]_-{\widetilde{f}}& \Sigma \ar@{=}[d]\\
  L' \ar@{->}[r]_-{\mu'} & \Sigma}
\] 
such that if $(\alpha,u_1,\dots,u_n,\beta)$ is a transition, then
$(f_0(\alpha),\widetilde{f}(u_1),\dots,\widetilde{f}(u_n),f_0(\beta))$
is a transition. The corresponding category is denoted by $\wts$.  The
$n$-transition $(\alpha,u_1,\dots,u_n,\beta)$ is also called a {\rm
  transition from $\alpha$ to $\beta$}.  The maps $f_0$ and
$\widetilde{f}$ will be also denoted by $f$. \ed

Every set $X$ may be identified with the weak transition system having
the set of states $X$, with no actions and no transitions.

It is usual in computer science to work in the comma category
$\{\iota\} \ddownarrow \wts$ where the image of the state $\iota$
represents the initial state of the process which is modeled. It then
makes sense to restrict to the states which are reachable from this
initial state by a path of transitions. Hence we introduce the
following definitions:

\bd Let $X$ be a weak transition system and let $\iota$ be a state of
$X$. A state $\alpha$ of $X$ is {\rm reachable from $\iota$} if it is
equal to $\iota$ or if there exists a finite sequence of transitions
$t_i$ of $X$ from $\alpha_i$ to $\alpha_{i+1}$ for $0\leq i \leq n$
with $n\geq 0$, $\alpha_0 = \iota$ and $\alpha_{n+1} = \alpha$.  \ed

\bd A {\rm star-shaped weak transition system} is an object $\{\iota\}
\rightarrow X$ of the comma category $\{\iota\} \ddownarrow \wts$ such
that every state of the underlying weak transition system $X$ is
reachable from $\iota$. The full subcategory of $\{\iota\} \ddownarrow
\wts$ of star-shaped weak transition systems is denoted by
$\wts_{\bullet}$. \ed

\bp \label{full-coref} The category $\wts_{\bullet}$ is a full
isomorphism-closed coreflective subcategory of $\{\iota\} \ddownarrow \wts$.  \ep

\bpf Let $\{\iota\} \to X$ be an object of $\{\iota\} \ddownarrow
\wts$. Let $\underline{T}_\iota(X)$ be the set of transitions
$(\alpha,u_1,\dots,u_n,\beta)$ of $X$ such that the initial state
$\alpha$ is reachable from $\iota$. Note that this implies that
$\beta$ is reachable from $\iota$ as well.  The set
$\underline{T}_\iota(X)$ satisfies the multiset axiom since permuting
the actions does not change the initial state of a transition. It also
satisfies the patching axiom because, with the notations of the
patching axiom in Definition~\ref{def_HDTS}, $\nu_1$ is reachable from
$\iota$. Therefore the triple consisting of the set of states of $X$
which are reachable from $\iota$, the set of actions of $X$ with the
same labelling map $\mu$, and the set of transitions
$\underline{T}_\iota(X)$ yields a well-defined weak transition system
$\underline{C}_\iota(X)$. By construction, the map $\iota \to
\underline{C}_\iota(X)$ is a star-shaped weak transition
system. Consider a commutative square
\[
\xymatrix
{
\{\iota\} \ar@{=}[rr]\fD{} && \{\iota\}\fD{} \\
&& \\
Y \ar@{->}[rr]^-{f} && X
}
\]
where $\{\iota\} \to Y$ is a star-shaped weak transition system. By
construction, for every state $\alpha$ of $Y$, $f(\alpha)$ is
reachable from $\iota$ and every transition
$(\alpha,u_1,\dots,u_n,\beta)$ of $Y$ is therefore mapped to a
transition $(f(\alpha),f(u_1),\dots,f(u_n),f(\beta))$ of
$\underline{T}_\iota(X)$. Therefore the commutative square above
factors uniquely as the composite of commutative squares
\[
\xymatrix {
  \{\iota\} \ar@{=}[rr]\fD{} && \{\iota\} \ar@{=}[rr]\fD{} &&  \{\iota\}\fD{} \\
  && \\
  Y \ar@{->}[rr]\ar@/_20pt/@{->}[rrrr]_-{f} &&
  \underline{C}_\iota(X)
  \fR{\subset} &&X.  }
\]
\epf

\bp \label{small-cone-inj} The category $\wts_{\bullet}$ is a small
cone-injectivity class of $\{\iota\} \ddownarrow \wts$ and all maps of the
cone can be chosen to be cofibrations. \ep

\bpf The category $\wts_{\bullet}$ is a small cone-injectivity class
with respect to the small cone formed by the inclusions of the weak
transition system $\{\iota,\alpha\}$ in the weak transition systems
\[\iota \stackrel{t_1} \longrightarrow \bullet \longrightarrow \dots
\longrightarrow \bullet \stackrel{t_n} \longrightarrow \alpha\] for
all $n\geq 0$ and all transitions $t_1,\dots,t_n$ with the labelling
map $\id_\Sigma$ (note that we must include in the cone the set map
$\{\iota,\alpha\}\to \{\iota\}$). The cone is small because there is a
\emph{set} of labels $\Sigma$. Finally, all maps of the cone are
cofibrations of weak transition systems because on the sets of
actions, they are all of them the inclusion of the empty set in some
set.  \epf

\begin{cor} \label{wts_bullet_locpre} The category $\wts_{\bullet}$ is locally
  presentable. \end{cor}

\bpf The category $\wts_{\bullet}$ is accessible by
\cite[Proposition~4.16]{MR95j:18001}. It is cocomplete by
Proposition~\ref{full-coref}. Therefore it is locally presentable.
\epf

We can now conclude this paper with the following application:

\bth There exists a left determined model structure on the category
$\wts_{\bullet}$ of star-shaped weak transition systems with respect
to the class of maps such that the underlying maps of weak transition
systems are one-to-one on actions.  \eth

\bpf[Sketch of proof] The category $\wts_{\bullet}$ is bicomplete by
Corollary~\ref{wts_bullet_locpre}.  By
\cite[Theorem~5.11]{cubicalhdts}, there exists an Olschok model
structure on the category of weak transition systems such that the
cofibrations are the maps which are one-to-one on actions. Let
$\cyl:\wts\to \wts$ be the cylinder functor which is described in the
proof of \cite[Proposition~5.8]{cubicalhdts}. The map $\{\iota\}
\sqcup \{\iota\} \to \cyl(\{\iota\})$ is epic because by
\cite[Proposition~5.8]{cubicalhdts}, one has $\cyl(\{\iota\}) =
\{\iota\}$. Hence we can apply Theorem~\ref{cartesian-comma}. We
obtain an Olschok model category on the comma category $\{\iota\}
\ddownarrow \wts$.  By Lemma~\ref{l0}, the cylinder of $\{\iota\}
\ddownarrow \wts$ is obtained by identifying two states in
$\cyl(X)$. Since the set of states of $\cyl(X)$ is equal to the set of
states of $X$ by the calculations made in the proof of
\cite[Proposition~5.8]{cubicalhdts}, the two states identified are
actually equal. This means that the underlying weak transition system
of $\cyl_{\{\iota\}}(\{\iota\} \to X)$ is $\cyl(X)$.  Therefore by
Corollary~\ref{far-fetched}, and since $\cyl$ is very good by
\cite[Proposition~5.7]{cubicalhdts}, the cylinder $\cyl_{\{\iota\}}$
is very good and the Olschok model category $\{\iota\} \ddownarrow
\wts$ is left determined. Let $\{\iota\}\to X$ be a star-shaped weak
transition system.  By the calculations made in the proof of
\cite[Proposition~5.8]{cubicalhdts} again, the set of actions of
$\cyl(X)$ is $L\p \{0,1\}$ where $L$ is the set of actions of $X$ and
a tuple $(\alpha,(u_1,\epsilon_1), \dots,(u_n,\epsilon_n),\beta)$ is a
transition of $\cyl(X)$ if and only if $(\alpha,u_1, \dots,u_n,\beta)$
is a transition of $X$. Therefore a state of $\cyl(X)$ is reachable
from $\{\iota\}$ if and only if it is reachable from $\{\iota\}$ in
$X$ (choose $\epsilon_i=0$ for all intermediate transitions).  One
deduces that if $\{\iota\} \to X$ is a star-shaped weak transition
system, then $\cyl_{\{\iota\}}(\{\iota\} \to X)$ is a star-shaped weak
transition system as well.  Using Proposition~\ref{small-cone-inj}, we
can now apply Theorem~\ref{cartesian-coreflective1}: we have obtained
an Olschok model structure which is left determined.  The proof is
complete.  \epf

\end{document}